\documentclass[12pt, final]{l4dc2023}

\title[Regret for Risk-aware Linear Quadratic Control]{Regret Analysis for Risk-aware Linear Quadratic Control}
\usepackage{times}
\usepackage{graphicx}
\usepackage{amsmath}
\usepackage{optidef}
\usepackage{amssymb}
\usepackage{multicol}
\usepackage{cuted}
\usepackage{multirow}
\usepackage{caption}
\usepackage{xcolor}
\usepackage{hyperref}
\usepackage{algorithm}
\usepackage{algpseudocode}
\usepackage{scalerel}
\usepackage{bm}
\definecolor{steelblue}{RGB}{70,130,180}

\usepackage{tikz}

\def\bbr{\mathbb R}

\def\bbs{\mathbb S}
\def\bbp{\mathbb P}
\def\bbe{\mathbb E}

\def\calp{\mathcal P}

\newcommand{\PP}{\mathbb{P}}

\newcommand{\norm}[1]{\left\|#1 \right\|}

\renewcommand{\norm}[1]{\left\lVert {#1} \right\rVert}

\DeclareMathOperator*{\argmin}{argmin} 
 
\newtheorem{assumption}{Assumption}

\author{%
 \Name{Venkatraman Renganathan} \Email{venkatraman.renganathan@control.lth.se}
 \AND
 \Name{Dongjun Wu} \Email{dongjun.wu@control.lth.se}\\
 \addr Department of Automatic Control - LTH\\
Lund University, Box 118, SE-221 00, Sweden.%
}

\begin{document}

\maketitle

\begin{abstract}%
This paper investigates the regret associated with the Distributionally Robust Control (DRC) strategies used to address multistage optimization problems where the involved probability distributions are not known exactly, but rather are assumed to belong to specified ambiguity families. We quantify the price (distributional regret) that one ends up paying for not knowing the exact probability distribution of the stochastic system uncertainty while aiming to control it using the DRC strategies. The conservatism of the DRC strategies for being robust to worst-case uncertainty distribution in the considered ambiguity set comes at the price of lack of knowledge about the true distribution in the set. We use the worst case Conditional Value-at-Risk to define the distributional regret and the regret bound was found to be increasing with tighter risk level. The motive of this paper is to promote the design of new control algorithms aiming to minimize the distributional regret.%
\end{abstract}

\begin{keywords}%
Regret Analysis, Distributionally Robust Control, Moment-based Ambiguity Sets%
\end{keywords}

\section{Introduction}
Stochastic optimization has been the backbone of success of several portfolio management systems in finance as they are equipped with necessary tools to handle risks of various kinds. Inspired by their application in finance, it is now being increasingly used in the control community as well. However, the success of stochastic optimization techniques depend on the knowledge of the true distribution of the random variable of interest as in \cite{woz:12}. For instance, it is a common practice in the well celebrated linear stochastic control approach to assume system uncertainties to be Gaussian as described in \cite{astrom:70, athans2013optimal}, just for the sake of simplicity and tractability. However, it is well known that system uncertainties in reality might be not known exactly and further they might be non-Gaussian. To address the shortcoming of not knowing the true distribution, the Distributionally Robust Optimization (DRO) approaches were developed in \cite{wiesemann2014distributionally, gao:22} to hedge against the worst case distributions of the uncertainties. In connection to the DRO approach, authors in \cite{kishida2022risk} proposed a risk aware linear quadratic control by hedging against the worst case distribution of the process disturbance using Conditional Value-at-Risk $\mathrm{(CVaR)}$. \\

In the process of being robust against any worst case arbitrary distribution in the assumed ambiguity set of distributions governing the random variable, all DRO techniques end up resulting in conservative decisions due to their pessimistic approach to the decision making problem. For instance, risk bounded motion planning for nonlinear robotic systems can be effectively performed as described in \cite{RENGANATHAN2023103812} if we exactly know how much unnecessary conservatism is being exerted by DRO approach while deciding the control actions. On the contrary, if we analyze the same problem through the optimistic lens, it naturally leads us to investigate the price that one would pay for rather hedging against the worst case distribution instead of the true distribution of the system uncertainties (referred later as \emph{distributional regret}). That is, we are interested in computing how much would any distributionally robust optimal control algorithm pay for not knowing the true distribution of the system uncertainty while trying to control the system by hedging against the worst case distribution of the system uncertainty. Out of several available DRO parameterizations, we address the regret with respect to moment based ambiguity set based problem formulation, especially with the first two moments as described in \cite{delage2010distributionally, jiang2018risk, nie2021distributionally, zymler2013distributionally}. The regret associated with other DRO parameterizations are equally interesting and open as well and are currently being pursued in our ongoing future work. \\

To compute the regret associated with the problem formulation with moment based ambiguity sets, it is natural to consider the DRO objective function formulated using the expectation or the worst case $\mathrm{CVaR}$ as the risk functional as described in \cite{rockafellar2000optimization, chapman2022risk}. However, we show that expectation is not the right risk functional to be considered in our problem setting as it does not effectively capture the tail probabilities from different distributions in the considered ambiguity set. Analysing the regret for not knowing the exact distributions has been carried out in the field of auction based game theory in \cite{guo2022auctions} where the authors used the learning approach to minimize the regret for not knowing the true distribution of items to buy. Our work is inspired by their motivation and we apply the same reasoning to the field of control theory for systems with stochastic uncertainties. Specifically, we study the regret incurred by the risk-aware linear quadratic control algorithm explained in \cite{kishida2022risk}. Similar analysis can be done in future for the risk-averse MPC described in \cite{sopasakis2019risk}.\\

\emph{Contributions:}  To the best of our knowledge, this is one of the first articles in the literature to analyze the regret of a control policy that rather hedges against the worst case distributions of the stochastic system uncertainties and not knowing the exact true distribution of stochastic system uncertainties. Our main contributions in this article are as follows:
\begin{enumerate}
    \item We define the distributional regret for discrete time stochastic linear systems with system uncertainties namely the initial state $x_0$ and the disturbance $w_{k}$ modelled using moment based ambiguity sets containing distributions consistent with the set of finite first two moments.
    \item Our distributional regret definition proposes the use of distributionally robust risk functional namely the worst-case Conditional Value-at-Risk $\mathrm{(CVaR)}$ to account for the price that a distributionally robust control algorithm would pay for not knowing the true distribution of the stochastic system uncertainties and purposefully hedging the worst case distribution instead.
    \item Our proposed regret analysis is simulated using numerical examples to demonstrate the growth of distributional regret over tighter risk level and time.
\end{enumerate}
Following a short summary of notations and preliminaries, the rest of the paper is organized as follows: We propose a regret definition for distributionally robust control with moment based ambiguity sets in \S \ref{sec_prob_formulation}. Subsequently, a regret analysis is performed in \S \ref{sec_dr_regret_analysis}. The simulation results demonstrating our proposed approach is explained in \S \ref{sec_num_sim}. Finally, the paper is closed in \S \ref{sec_conclusion} along with a summary and directions for future research.


\section*{Notations \& Preliminaries}
The set of real numbers and natural numbers are denoted by $\mathbb{R}$ and $\mathbb{N}$, respectively. For $x \in \mathbb{R}$, we denote by $x^{+} := \max\{x, 0\}$. An identity matrix in dimension $n$ is denoted by $I_{n}$. For a non-zero vector $x \in \bbr^{n}$ and a matrix $P \in \bbs^{n}_{++}$ or $P \succ 0$, we denote $\left \| x \right \Vert^{2}_{P} = x^{\top} P x$, where $\bbs^{n}_{++}$ is the set of all positive definite matrices. We say $P \succ Q$ or $P \succeq Q$ if $P - Q \succ 0$ or $P - Q \succeq 0$ respectively. For a vector $x \in \mathbb{R}^{n}$ and a matrix $Q \in \mathbb{R}^{n \times n}$, the notation $\norm{x}_{Q}$ denotes its scaled quadratic form $\sqrt{x^{\top} Q x}$. 
We denote by $\mathcal{B}(\bbr^{d})$ and $\mathcal{P}(\bbr^{d})$ the Borel $\sigma-$algebra on $\bbr^{d}$ and the space of probability measures on $(\bbr^{d}, \mathcal{B}(\bbr^{d}))$ respectively. A probability distribution with mean $\mu$ and covariance $\Sigma$ is denoted by $\mathbb{P}(\mu, \Sigma)$, and specifically $\mathcal{N}_{d}(\mu, \Sigma)$ if the distribution is normal in $\mathbb{R}^{d}$. 
Given a constant $q \in \bbr_{\geq 1}$, the set of probability measures in $\mathcal{P}(\bbr^{d})$ with finite $q-$th moment is denoted by $\mathcal{P}_{q}(\bbr^{d}) := \left\{ \mu \in \mathcal{P}(\bbr^{d}) \mid \int_{\bbr^{d}} \left \Vert x \right \|^{q} d\mu < \infty \right\}$. 


\section{Control of Linear Systems with Unknown Distributions for System Uncertainties} \label{sec_prob_formulation}
Consider a linear, stochastic, discrete time-invariant system as follows 
\begin{equation} \label{eqn_dt_stochastic_sys}
x_{k+1}= A x_{k} + B u_{k} + D w_{k}, \quad k = 0,\dots,T-1,
\end{equation}
where $x_{k} \in \bbr^n$ and $u_{k} \in \bbr^m$ is the system state and input at time $k$, respectively and $T$ denotes the total time horizon. Further, $w_{k} \in \bbr^r$ denotes the stochastic process noise. Further, $A \in \mathbb{R}^{n \times n}$ is the dynamics matrix, $B \in \mathbb{R}^{n \times m}$ is the input matrix and $D \in \mathbb{R}^{n \times r}$ is the disturbance matrix and the pair $(A,B)$ is assumed to be stabilizable. We denote by $\mathbf{w} := ( w_{0}, \dots,w_{T-1})$.

\subsection{Problem Formulation Using Moment-based Ambiguity Set} 
Based on the assumption made on the knowledge of the distributions of the stochastic system uncertainties $x_0$ and $w$, we will define and analyse regret associated with it. There are several ambiguity set definitions for modeling the stochastic system uncertainties such as moment based, Wasserstein metric based, multimodality of the distribution and other parameterization. However, we limit our analysis to only moment-based ambiguity set. The analysis of regret for other parameterizations are left as an interesting direction of research for the future. The preliminaries of the details of the ambiguity sets used in this paper are given below.
\label{subsec_MBAS}
\begin{assumption} \label{assume_1}
The distribution of the disturbance $w_{k}$ is the same $\forall k = 0,\dots, T-1$ meaning that $w_{k}$ is identically distributed across time and so $\mathbb{P}_{w_{k}} = \mathbb{P}_{w}$. Further, $\mathbb{E}[w_{k}w^{\top}_{j}] = 0, \forall k \neq j$.    
\end{assumption}
The distribution of $w_{k}$, namely $\mathbb{P}_{w}$, is unknown but is believed to be belonging to a moment-based ambiguity set of distributions, $\calp^w$ given by
\begin{equation} \label{eqn_ambig_w}
\calp^w = \left\{ \bbp_{w} \in \mathcal{P}_{2}(\bbr^{r}) \mid \bbe[w] = 0, \, \bbe[w w^{\top}] = \Sigma_{w} \right\}.
\end{equation}
Note that there are infinitely many distributions present in the considered set $\calp^w$. For instance, both multivariate Gaussian and multivariate Laplacian distributions with zero mean and covariance $\Sigma_{w}$ belong to $\calp^w$ with the latter having heavier tails than the former. We assume that the system is controllable under zero process noise meaning that given any $x_{0}, x_{f} \in \mathbb{R}^{n}$, there exists a sequence of control input $\{u_{k}\}^{N-1}_{k=0}$ that steers the system state from $x_{0}$ to $x_{f}$ for large $N \in \mathbb{N}$. Since the considered system is stochastic, the initial condition $x_{0}$ is subject to a similar uncertainty model as the process noise, with the distribution belonging to a moment-based ambiguity set, $\calp^{x_0}$ given by 
\begin{equation} \label{eqn_ambig_x_0}
    \calp^{x_0} = \left\{ \PP_{x_{0}} \in \mathcal{P}_{2}(\bbr^{n}) \mid \bbe[x_0] = \mu_{x_{0}}, \bbe[(x_0-\mu_0) (x_0-\mu_0)^{\top}] = \Sigma_{x_{0}} \right\}.
\end{equation}
If the control law $u_{k}$ is selected as an affine function of the state $x_{k}$ at any time $k$, similar moment-based ambiguity sets $\calp^{x_{k}}$ can be written for $\mathbb{P}_{x_{k}}, \forall k \in 1,\dots,N$ using the propagated mean and covariance at time $k$. Note that $\calp^{x_{k}}$ would \emph{not} be empty $\forall k \in [1,N]$ as Gaussian distribution (due to its invariance under linear transformation) would be a guaranteed member in that set. Further, $x_0$ is assumed to independent of $w_{k}, \forall k$. Here, we present the regret analysis for the control problem involving moment-based ambiguity set-based problem formulation. Towards the regret definition, we require certain assumptions on the control policy to be used by the system considered in \eqref{eqn_dt_stochastic_sys}.
\begin{assumption}
We restrict the class of control policies to be causal and linear state feedback policies so that at each time instant $k$, the control input $u_{k} = \pi_{k}(x_{k}, x_{k-1}, \dots, x_{0})$, where the policy $\pi_{k} \in \Pi$ and $\Pi$ denotes the set of all linear state feedback policies.
\end{assumption}

\subsubsection{Minimax Problem with Moment-based Ambiguity Set}
For brevity of notation, let us denote by $\mathbf{x} := (x_{0}, x_{1}, \dots, x_{T})$ and $\mathbb{P}_{\mathbf{x}} := (\mathbb{P}_{x_{0}}, \mathbb{P}_{x_{1}}, \dots, \mathbb{P}_{x_{T}})$. Consider, the following \emph{random} quadratic cost function
\begin{equation} \label{eqn_moment_cost}
    J_{T}\left(\pi, \mathbf{x}, \mathbf{w}\right) = 
    \norm{x_{T}}^{2}_{Q}
    + \sum_{k=0}^{T-1} \left(\norm{x_{k}}^{2}_{Q} + \norm{u_{k}}^{2}_{R} \right),
\end{equation}
where $T \in \mathbb{N}$ is the finite time horizon, $Q \succeq 0$ and $R \succ 0$. That is, \eqref{eqn_moment_cost} is the resulting cost when the state evolution of the system given by \eqref{eqn_dt_stochastic_sys} starting from $x_{0} \sim \mathbb{P}_{x_{0}} \in \mathcal{P}^{x_{0}}$ is governed by the policy $\pi \in \Pi$ under the disturbances $w_{k}, \forall k = 0,\dots,T-1$ coming from the distribution $\mathbb{P}_{w} \in \mathcal{P}^{w}$ and the resulting distribution of states being $\mathbb{P}_{\mathbf{x}}$. We define the risk functional $\rho$ operating on the random cost \eqref{eqn_moment_cost} to be a map from the space of random variables to $\mathbb{R} \cup \{\infty\}$. There are several risk functionals available from the finance literature described in \cite{rockafellar2000optimization} such as the expectation, Value-at-Risk, Conditional Value-at-Risk, mean-variance etc. A good risk functional for a specific application should satisfy certain required axioms mentioned in \cite{majumdarrisk}. In this problem setting, risk functionals such as the expectation, Value-at-Risk, Conditional value-at-Risk can be used. Consider a random variable $Z \in \mathbb{R}$ with finite second order moments modeled using moment based ambiguity set $Z \sim \mathbb{P}_{Z} \in \mathcal{P}^{Z} := \{ \mathbb{P}_{Z} \mid \mathbb{E}[Z] = \mu_{Z}, \mathbb{E}[(Z-\mu_{Z})^{2}] = \Sigma_{Z} \}$. Then, $\mathrm{VaR}$ quantifies the threshold for a given tail probability corresponding to the distribution of the random variable and $\mathrm{CVaR}$ quantifies the conditional expectation of losses exceeding the threshold defined using $\mathrm{VaR}$. That is, for a given risk level $\alpha \in (0,1)$, we define
\begin{align}
\mathrm{VaR}_{\alpha}(Z) 
&:= \inf\{z \in \mathbb{R} \mid \mathbb{F}(z) \geq 1 - \alpha\}, \quad \text{and} \label{eqn_var},\\
\mathrm{CVaR}_{\alpha}(Z) 
&:= \inf_{s \in \mathbb{R}} \left( s+\frac{1}{\alpha} \mathbb{E}[(Z - s)^{+}] \right)
= \mathrm{VaR}_{\alpha}(Z) + \frac{1}{\alpha} \mathbb{E}\left[(Z - \mathrm{VaR}_{\alpha}(Z))^{+}\right],
\label{eqn_cvar}
\end{align}
where \eqref{eqn_cvar} is convex in $Z$ and $\mathbb{F}(z)$ denotes the cumulative distribution function of the random variable $Z$ and $\mathbb{F}(z)$ is both right-continuous and non-decreasing function. One can find the control policy that minimizes \eqref{eqn_moment_cost}, if the distribution of the system uncertainties namely $\mathbb{P}_{\mathbf{x}}, \mathbb{P}_{w}$ are known exactly as it would aide us in computing the risk functionals. However, in reality the distribution of the random variable is not known exactly or only partial information is known. Then in that case, we define the worst case $\mathrm{CVaR}_{\alpha}(Z)$ of the random variable $Z$ as 
\begin{align}
\overline{\mathrm{CVaR}}_{\alpha}(Z) 
&:= \sup_{\mathbb{P}_{Z} \in \mathcal{P}^{Z}} \mathrm{CVaR}_{\alpha}(Z) = \inf_{s \in \mathbb{R}} \left( s+\frac{1}{\alpha} \sup_{\mathbb{P}_{Z} \in \mathcal{P}^{Z}} \mathbb{E}[(Z - s)^{+}] \right),
\label{eqn_worst_cvar}
\end{align}
where we interchanged the $\inf$ and $\sup$ operation using the stochastic saddle point theorem in \cite{shapiro2002minimax}. Note that when $\alpha \rightarrow 1$, $\overline{\mathrm{CVaR}}_{\alpha}(Z)$ reduces to $\mathbb{E}[Z]$. For instance, the classical stochastic linear quadratic regulator (LQR) problem setting assumes that the disturbances $w_{k} \sim \mathcal{N}(0, \Sigma_{w}), \forall k$ with $\Sigma_w \in \bbs^{r}_{+}$ and the initial state $x_0 \sim \mathcal{N}(\mu_0, \Sigma_{0})$ with covariance $\Sigma_0 \in \bbs^{n}_{+}$. Due to the linearity of the system dynamics \eqref{eqn_dt_stochastic_sys}, it is straightforward to see that $\mathbb{P}_{x_{k}}, \forall k = 0,1,\dots,T$ will be Gaussian as well.
Then, LQR amounts to solving the following optimization problem under the full state feedback control law given by $u_{k} = \pi_{k} (x_{k})$: 
\begin{subequations} \label{eqn_stochasticLqr}
\begin{alignat}{2}
\underset{\pi \in \Pi}{\textrm{minimize}} \quad & \mathbb{E}\left[J_{T}\left(\pi, \mathbf{x}, \mathbf{w}\right) \right]\\
\textrm{subject to} \quad & x_{k+1}= A x_{k} + B u_{k} + D w_{k}, \quad k = 0,\dots,T-1,\\
&x_{0} \sim \mathcal{N}(\mu_0, \Sigma_{0}), w_{k} \sim \mathcal{N}(0, \Sigma_{w}), \forall k = 0,\dots,T-1.
\end{alignat}
\end{subequations}
\noindent The optimal trajectory of the system \eqref{eqn_dt_stochastic_sys} denoted by $\bar{\mathbf{x}} =(\bar{x}_0, \bar{x}_1, \cdots, \bar{x}_{T})$ satisfying \eqref{eqn_stochasticLqr} is obtained under the optimal control input $\bar{u}_{k} = \bar{\pi}_{k}(\bar{x}_{k}) = \bar{K}_{k} \bar{x}_{k}$, with $\bar{K}_{k} \in \mathbb{R}^{m \times n}$ and the policy $\bar{\pi} \in \Pi$ that minimizes \eqref{eqn_stochasticLqr} can be found using standard dynamic programming techniques which can be found in \cite{bertsekasbook}. Now consider the variation of stochastic LQR problem where the system uncertainties $w$ and $x_0$ are modelled using the moment-based ambiguity sets \eqref{eqn_ambig_w} and \eqref{eqn_ambig_x_0} respectively as described in \cite{kishida2022risk}:
\begin{subequations} \label{eqn_momentbasedLqr}
\begin{alignat}{2}
\underset{\pi \in \Pi}{\textrm{minimize}} \quad & \underset{\mathbb{P}_{\mathbf{x}}, \mathbb{P}_{w}}{\sup} \rho\left[J_{T}\left(\pi, \mathbf{x}, \mathbf{w}\right) \right]\\
\textrm{subject to} \quad & x_{k+1}= A x_{k} + B u_{k} + D w_{k}, \quad k = 0,\dots,T-1,\\
&x_{0} \sim \mathbb{P}_{x_{0}} \in \mathcal{P}^{x_{0}}, w_{k} \sim \mathbb{P}_{w} \in \mathcal{P}^{w}, \forall k = 0,\dots,T-1,
\end{alignat}
\end{subequations}
\noindent where the risk functional $\rho$ is \emph{not} the expectation operator. Let the optimal feedback control input at time $k$ be $u^{\star}_{k} = \pi^{\star}_{k}(x^{\star}_{k}) = K^{\star}_{k} x^{\star}_{k}$, with $K^{\star}_{k} \in \mathbb{R}^{m \times n}$ given by the policy $\pi^{\star}_{k} \in \Pi$\footnote{Note that we restrict ourselves to set of linear state feedback policies though potentially a nonlinear state feedback policy can do better that the linear counterparts.}. Then, the optimal trajectory of the system \eqref{eqn_dt_stochastic_sys} is denoted by $\mathbf{x}^{\star} = (x_0^{\star}, x_1^{\star}, \cdots, x_{T}^{\star})$ and is achieved under the policy $\pi^{\star} \in \Pi$ given by 
\begin{equation} \label{eqn_pi_star}
    \pi^{\star} = \argmin_{\pi \in \Pi} \sup_{ \mathbb{P}_{\mathbf{x}}, \mathbb{P}_{w}} \rho\left[ J_{T}\left(\pi, \mathbf{x}, \mathbf{w}\right) \right].
\end{equation}



\section{Regret of Distributionally Robust Control with Moment Based Ambiguity Sets} \label{sec_dr_regret_analysis}
We denote by $\mathbf{w}^{\star} = \{w^{\star}_{0},\dots,w^{\star}_{T-1}\}$ and $\breve{\mathbf{w}} = \{\breve{w}_{0},\dots,\breve{w}_{T-1}\}$ the worst case disturbances and the true disturbances realized from $\mathbb{P}^{\star}_{w} \in \mathcal{P}^{w}$ and $\breve{\mathbb{P}}_{w} \in \mathcal{P}^{w}$ respectively. Similarly, $\breve{\mathbf{x}} := \{ \breve{x}_{0}, \dots, \breve{x}_{T} \}$ denotes the true states with uncertainties stemming from their respective true distributions. Regret for online controllers can be addressed using \emph{policy regret} approach where the cost incurred by the online control policy operating without knowing the quantity of interest is compared against that of the optimal control policy available in hindsight that operates by knowing the quantity of interest. Since the distribution of the stochastic system uncertainties in stochastic LQR are exactly known, and moreover, it being the optimal linear state feedback solution, there is nothing that stochastic LQR can do to improve more as it has the entire knowledge about the stochastic system uncertainties. On the other hand, the case with the moment based ambiguity set is not straightforward and they demand a certain level of conservatism to be robust against any distribution in the respective ambiguity set. Hence, to investigate the conservatism of such control policies from an optimistic lens, we study their regret. Since, the policy $\pi^{\star}_{k} \in \Pi$ minimizes \eqref{eqn_momentbasedLqr}, it is important to mention in what sense it is optimal. That is, the policy $\pi^{\star} \in \Pi$ might not know the true distributions of the states and disturbance namely $\breve{\mathbb{P}}_{x_{k}} \in \mathcal{P}^{x_{k}}, \forall k = 0,\dots, T$, and $\breve{\mathbb{P}}_{w}$ respectively but it might hedge against the worst case distributions of the states and disturbance namely $\mathbb{P}^{\star}_{x_{k}} \in \mathcal{P}^{x_{k}}, \forall k = 0,\dots, T$, and $\mathbb{P}^{\star}_{w}$ respectively. Note that at any time $k = 0, \dots, T$, the true distribution of the state namely $\breve{\mathbb{P}}_{x_{k}} \in \mathcal{P}^{x_{k}}$ need not be equal to the worst case distribution of the state $\mathbb{P}^{\star}_{x_{k}} \in \mathcal{P}^{x_{k}}$. If $\pi^{\star} \in \Pi$ was designed to be robust against the worst case distributions of the states and the disturbance, then it will not fair well against the true distributions of the states and the disturbance, meaning that $\pi^{\star} \in \Pi$ will incur more cost for the system in face of the true distributions of the states and the disturbance which are not necessarily equal to the worst case distribution of the states and the disturbance. To quantify the price that the policy $\pi^{\star} \in \Pi$ will incur in such situations, one would expect to define the pseudo distributional regret as follows:
\begin{align} \label{eqn_regret_moment_ambig_almost}
    \bar{\mathbf{R}}(T, \pi^{\star}) := \rho\left[J_{T}\left(\pi^{\star}, \breve{\mathbf{x}}, \breve{\mathbf{w}}\right)\right] 
    - \rho\left[J_{T}\left(\pi^{\star}, \mathbf{x}^{\star}, \mathbf{w}^{\star}\right) \right],
\end{align}
where the risk functional $\rho$ is taken to be the expectation operator $\mathbb{E}$ and the notations $\breve{\mathbb{P}}_{\mathbf{x}} := (\breve{\mathbb{P}}_{x_{0}}, \breve{\mathbb{P}}_{x_{1}}, \dots, \breve{\mathbb{P}}_{x_{T}})$ and $\mathbb{P}^{\star}_{\mathbf{x}} := (\mathbb{P}^{\star}_{x_{0}}, \mathbb{P}^{\star}_{x_{1}}, \dots, \mathbb{P}^{\star}_{x_{T}})$. Here, $\bar{\mathbf{R}}(T, \pi^{\star})$ denotes the pseudo distributional regret defined by the system's running cost incurred by the policy $\pi^{\star} \in \Pi$ for rather hedging against the worst case distributions and not being prepared against the true distributions of the stochastic system uncertainties.

\begin{lemma}
The pseudo distributional regret given by \eqref{eqn_regret_moment_ambig_almost} experienced by the stochastic LQR policy $\bar{\pi} \in \Pi$ is zero. This means that there is nothing that stochastic LQR can do to improve more as it has the entire knowledge about the stochastic system uncertainties.  
\end{lemma}
\begin{proof}
Recall that in the stochastic LQR setting, $\mathbb{P}^{\star}_{\mathbf{x}} = \breve{\mathbb{P}}_{\mathbf{x}}$ and they are known to be Gaussian distributed. Similarly, $\mathbb{P}^{\star}_{w} = \breve{\mathbb{P}}_{w}$ and they are known to be Gaussian distributed. Substituting these inferences into the regret definition given by \eqref{eqn_regret_moment_ambig_almost} with $\bar{\pi}$ replacing $\pi^{\star}$, we immediately see that $\bar{\mathbf{R}}(T, \bar{\pi}) = 0$. 
\end{proof}
A more careful observation will depict that using the expectation as a risk functional acting on the random cost to evaluate the risk of using the policy $\pi^{\star} \in \Pi$ when not knowing the true distribution is not a good idea given the moment based ambiguity set modeling approach. The following lemma describes the observation precisely. 
\begin{lemma}
The pseudo distributional regret defined using \eqref{eqn_regret_moment_ambig_almost} by policy $\pi^{\star} \in \Pi$ is zero even when the true distribution of the states $\breve{\mathbb{P}}_{x_{k}} \in \mathcal{P}^{x_{k}}, \forall k = 0,\dots,T$ and true distribution of the disturbance $\breve{\mathbb{P}}_{w} \in \mathcal{P}^{w}$ are different from their worst case counterparts. 
\end{lemma}
\begin{proof}
Given that policy $\pi^{\star} \in \Pi$ is linear in system states, the control input at any time $k$ is given by $u^{\star}_{k} = K^{\star}_{k} x^{\star}_{k}$. Then, expanding the definition of $\bar{\mathbf{R}}(T, \bar{\pi})$ using the control input $u^{\star}_{k}$ and the fact that $\mathbb{E}[x^{\top} Q x] = \mathbf{tr}(Q \Sigma_{x}) + \mu^{\top}_{x} Q \mu_{x}$, we get 
\begin{align*}
\bar{\mathbf{R}}(T, \bar{\pi}) 
&= \left[ \sum^{T}_{k = 0} \left( \mathbf{tr}(Q^{\star}_{k} \breve{\Sigma}_{x_{k}}) +
\norm{\mathbb{E}[\breve{x}_{k}]}^{2}_{Q^{\star}_{k}}
\right)\right] - \left[ \sum^{T}_{k = 0} \left( \mathbf{tr}(Q^{\star}_{k} \Sigma^{\star}_{x_{k}}) + 
\norm{\mathbb{E}[x^{\star}_{k}]}^{2}_{Q^{\star}_{k}}
\right)\right] = 0,
\end{align*}
where $Q^{\star}_{k} = Q + {K^{\star}_{k}}^{\top} R K^{\star}_{k}$ for $k = 0,\dots, T-1$ and $Q^{\star}_{T} = Q$ and further we used the fact that for $\{\mathbb{P}^{\star}_{x_{k}}, \breve{\mathbb{P}}_{x_{k}}\}^{T}_{k=0} \in \mathcal{P}^{x_{k}}$, the associated mean and covariances are the same due to moment based ambiguity set formulation.
\end{proof}

\subsection{Definition of Distributional Regret}
This shows that the regret defined using \eqref{eqn_regret_moment_ambig_almost} will not qualitatively capture the cost for not knowing the true distributions of the stochastic uncertainties. This forces us to look for another risk functional to evaluate the policy $\pi^{\star} \in \Pi$ when the true distributions of the stochastic system uncertainties are not known similar to the works of \cite{whittle1981risk}. In light of the considered moment based ambiguity sets, we propose to define the regret using distributionally robust risk functional such as the Value-at-Risk (VaR) and Conditional Value-at-Risk (CVaR). Recall that, both multivariate Gaussian and multivariate Laplacian distributions with zero mean and covariance $\Sigma_{w}$ belong to $\calp^w$ belong to the same ambiguity set but the latter having heavier tails than the former. This means that if the uncertainty is coming from a true distribution having a heavier tail, then expectation would not be a good metric to evaluate the risk and distributionally robust risk functional such as the VaR, CVaR are needed. While it is compelling to use the $\mathrm{VaR}_{\alpha}(\cdot)$ to define the regret as it is closely related to the chance constraints associated with the given moments, it is not a good choice due to $\mathrm{VaR}_{\alpha}(\cdot)$ being not a coherent risk measure as it does not satisfy all the required axioms (specifically sub-additivity) described in \cite{majumdarrisk} and this would eventually thwart our efforts to find any regret bound later. Further, we do not know the distribution of the random cost. Hence, we finally propose to use the worst-case $\mathrm{CVaR}_{\alpha}(\cdot)$ to define the distributional regret as it is a coherent risk functional as described in Proposition 1 of \cite{zhu2009worst}.

\begin{definition}
The distributional regret associated with the policy $\pi^{\star} \in \Pi$ defined using \eqref{eqn_pi_star} for a given risk level $\alpha \in (0,1)$ and for hedging against the worst case distribution and not knowing the true distributions of the stochastic system uncertainties is defined \footnote{We can describe the \emph{Knowledge Gap} incurred by the policy $\pi$ using the difference of actual uncertainty and the worst case uncertainty experienced by the system. That is, $\mathbf{KG}(\pi, T) = \sum^{T}_{k=0} f(\Delta w_{k}, \Delta x_{k})$ where, $\Delta w_{k} : = \breve{w}_{k} - w^{\star}_{k}$ and $\Delta x_{k} : = \breve{x}_{k} - x^{\star}_{k}$. Analysing regret due to this knowledge gap is an interesting direction for future research.} as
\begin{subequations}
\label{eqn_regret_moment_ambig}
\begin{align}
    \mathbf{R}_{\alpha}(T, \pi^{\star}) 
    &:= \sup_{ \mathbb{P}_{\mathbf{x}}, \mathbb{P}_{w}} \mathrm{CVaR}_{\alpha}\left[J_{T}\left(\pi^{\star}, \breve{\mathbf{x}}, \breve{\mathbf{w}}\right)\right] 
    - \sup_{ \mathbb{P}_{\mathbf{x}}, \mathbb{P}_{w}} \mathrm{CVaR}_{\alpha}\left[J_{T}\left(\pi^{\star}, \mathbf{x}^{\star}, \mathbf{w}^{\star}\right) \right] \\
    &:=
    \overline{\mathrm{CVaR}}_{\alpha}\left[J_{T}\left(\pi^{\star}, \breve{\mathbf{x}}, \breve{\mathbf{w}}\right)\right] 
    - \overline{\mathrm{CVaR}}_{\alpha}\left[J_{T}\left(\pi^{\star}, \mathbf{x}^{\star}, \mathbf{w}^{\star}\right) \right].
\end{align}
\end{subequations}
\end{definition}
Note that the risk functional $\rho$ is chosen to be the worst-case $\mathrm{CVaR}_{\alpha}(\cdot)$ in the above regret definition. Here, the distributional regret is defined by the system's running cost incurred by the policy $\pi^{\star} \in \Pi$ for not being prepared against the true distributions of the stochastic system uncertainties. This helps us to study the distributional regret $\mathbf{R}_{\alpha}(T, \pi^{\star}) $ for the policy $\pi^{\star} \in \Pi$ for cases when $\mathbb{P}^{\star}_{x_{k}} \neq \breve{\mathbb{P}}_{x_{k}}, \forall k = 0, \dots, T$ despite the associated moments in the corresponding ambiguity set $\mathcal{P}^{x_{k}}$ being the same. That is, when the system is started from two different initial states, one sampled from an assumed worst case distribution lying in $\mathcal{P}^{x_{0}}$ and another one sampled from the true (possibly heavier tail) distribution lying in $\mathcal{P}^{x_{0}}$, the system states due to the later case will evolve differently (such as having a heavier tail) despite having same moments like the former case and the risk functional worst-case $\mathrm{CVaR}_{\alpha}(\cdot)$ has the capability to capture the difference and hence the qualitative distributional regret stemming from such a setting.  

\subsection{Regret Analysis}
To explicitly study the difference of the evolution of the states under the worst-case distribution and true distribution, we define the state error as the difference between the states evolved under the true stochastic system uncertainties and the worst case stochastic system uncertainties. That is, $\forall k = 0, \dots, T, e_{k} := \breve{x}_{k} - x^{\star}_{k}$. Note that at any time $k = 0,\dots, T$, we see that $\mathbb{E}[e_{k}] = \mathbb{E}[\breve{x}_{k} - x^{\star}_{k}] 
= \mathbb{E}[\breve{x}_{k}] - \mathbb{E}[x^{\star}_{k}] 
= 0$. Let us denote the second order moment matrix of the $e_{k}$ at time $k$ as
$\Omega_{k} = \begin{bmatrix} \Sigma_{e_{k}} & 0 \\ 0 & 1 \end{bmatrix}$, where $\Sigma_{e_{k}} 
:= \mathbb{E}[(\breve{x}_{k} - x^{\star}_{k}) (\breve{x}_{k} - x^{\star}_{k})^{\top}] 
= 2 \left( \mathbb{E}[\breve{x}_{k} \breve{x}^{\top}_{k}] -  \mathbb{E}[x^{\star}_{k}] \mathbb{E}[\breve{x}^{\top}_{k}] \right) = 2 \Sigma_{x_{k}}$.
\begin{lemma} \label{lemma_cvar_of_quadratic_cost}
For a zero-mean random vector $e_{k} \in \mathbb{R}^{n}$ at time $k$ with covariance $\Sigma_{e_{k}} \in \mathbb{R}^{n \times n}$, and a quadratic cost function $f(e_{k}) = \norm{e_{k}}^{2}_{P}$ with $P \succeq 0$, we see that for risk level $\alpha \in (0,1)$,
\begin{align}
    \overline{\mathrm{CVaR}}_{\alpha} (\norm{e_{k}}^{2}_{P}) = \frac{1}{\alpha} \mathbf{Tr}(\Sigma_{e_{k}} P). 
\end{align}
\end{lemma}
\begin{proof}
Proof is the same as Corollary 3.3 in \cite{kishida2022risk} with $f(e_{k}) = \norm{e_{k}}^{2}_{P}$.    
\end{proof}
\begin{lemma} \label{lemma_worst_cvar_linear_term}
(Theorem 21 of \cite{zymler2013distributionally}) Consider $g(e_{k}) = e^{\top}_{k} P e_{k} + 2  q^{\top} e_{k} + r$ with $P \succeq 0$ for the random vector $e_{k} \in \mathbb{R}^{n}$. Then, $\overline{\mathrm{CVaR}}_{\alpha}(g(e_{k}))$ for risk level $\alpha \in (0,1)$ is given by
\end{lemma}
\begin{subequations} \label{eqn_worst_cvar_linear_term}
\begin{alignat}{2}
\underset{s}{\inf} \quad & s + \frac{1}{\alpha} \mathbf{Tr}(\Omega_{k} M)\\
\textrm{subject to} \quad & M - \begin{bmatrix} P & q \\ q^{\top} & r-s    
\end{bmatrix} \succeq 0,\\
& M \in \mathbb{S}^{n+1}_{+}, s \in \mathbb{R}.
\end{alignat}
\end{subequations}

\begin{theorem} \label{thm_main}
The upper bound on the distributional regret defined using the risk functional worst-case $\mathrm{CVaR}_{\alpha}(\cdot)$ with risk level $\alpha \in (0,1)$ in \eqref{eqn_regret_moment_ambig} for policy $\pi^{\star} \in \Pi$ grows with the time $T$ and a tighter risk level $\alpha$ when the true distribution of the states $\breve{\mathbb{P}}_{x_{k}} \in \mathcal{P}^{x_{k}}, \forall k = 0,\dots,T$ and true distribution of the disturbance $\breve{\mathbb{P}}_{w} \in \mathcal{P}^{w}$ are different from their respective worst-case counterparts. 
\end{theorem}
\begin{proof}
Note that policy $\pi^{\star} \in \Pi$ is linear in system states with feedback gain matrix being $K^{\star}_{k} \in \mathbb{R}^{m \times n}$ at any time $k$. Then, expanding \eqref{eqn_regret_moment_ambig} with $Q^{\star}_{k} = Q + {K^{\star}_{k}}^{\top} R K^{\star}_{k}$, for $k = 0,\dots, T-1$ and $Q^{\star}_{T} = Q$ and eliminating $\breve{x}_{k}$ by using $\breve{x}_{k} = x^{\star}_{k} + e_{k}$, we get 
\begin{align}
\mathbf{R}_{\alpha}(T, \pi^{\star})
= \, & 
\overline{\mathrm{CVaR}}_{\alpha}\left[ \sum_{k=0}^{T} \norm{\breve{x}_{k}}^{2}_{Q^{\star}_{k}}
\right] - \overline{\mathrm{CVaR}}_{\alpha}\left[ \sum_{k=0}^{T} \norm{x^{\star}_{k}}^{2}_{Q^{\star}_{k}} \right]\\
= \, & 
\overline{\mathrm{CVaR}}_{\alpha}\left[ \sum_{k=0}^{T} \norm{{x}_{k}^{\star}}^{2}_{Q^{\star}_{k}} + \norm{e_k}_{Q_k^\star}^2 + 2 (Q_k^\star x_k^\star)^\top e_k \right] -  \overline{\mathrm{CVaR}}_{\alpha}\left[ \sum_{k=0}^{T}\norm{x^{\star}_{k}}^{2}_{Q^{\star}_{k}} \right],  \\
\leq \, & 
\overline{\mathrm{CVaR}}_{\alpha}\left[  \sum_{k=0}^{T}\norm{{x}_{k}^{\star}}^{2}_{Q^{\star}_{k}} \right] + \sum_{k=0}^{T}\overline{\mathrm{CVaR}}_{\alpha}\left[ \norm{e_k}_{Q_k^\star}^2 \right] + \sum_{k=0}^{T}\overline{\mathrm{CVaR}}_{\alpha} \left[2 (Q_k^\star x_k^\star)^\top e_k \right] \nonumber \\
 &  \quad - \overline{\mathrm{CVaR}}_{\alpha}\left[ \sum_{k=0}^{T}\norm{x^{\star}_{k}}^{2}_{Q^{\star}_{k}} \right],  \\
= \, &
\sum_{k=0}^{T} \overline{\mathrm{CVaR}}_{\alpha} \left[ 
\norm{e_{k}}^{2}_{Q^{\star}_{k}} \right]  
+ \sum_{k=0}^{T} \overline{\mathrm{CVaR}}_{\alpha}\left[2 
(Q^{\star}_{k} x^{\star}_{k})^{\top} e_{k} \right], \\
= \, &
\sum_{k=0}^{T} \frac{1}{\alpha} \mathbf{Tr}(\Sigma_{e_{k}} Q^{\star}_{k}) + \sum_{k=0}^{T} 
\underbrace{\overline{\mathrm{CVaR}}_{\alpha}\left[2 (Q^{\star}_{k} x^{\star}_{k})^{\top} e_{k} \right]}_{= G_{k}}, \label{eqn_cvar_upper_bound_bigger}
\end{align}
\noindent where Lemma \ref{lemma_cvar_of_quadratic_cost} and the sub-additivity property of the risk functional $\overline{\mathrm{CVaR}}_{\alpha}(\cdot)$ were utilised to get the required inequality \eqref{eqn_cvar_upper_bound_bigger}. Note that the terms $G_{k}, \forall k = 0,\dots,T$ in \eqref{eqn_cvar_upper_bound_bigger} can be found by solving the optimization problem \eqref{eqn_worst_cvar_linear_term} given in Lemma \ref{lemma_worst_cvar_linear_term} by substituting $P = 0_{n \times n}, q = Q^{\star}_{k} x^{\star}_{k}$ and $r = 0$. Since, the true distribution of the states $\breve{\mathbb{P}}_{x_{k}} \in \mathcal{P}^{x_{k}}, \forall k = 0,\dots,T$ and true distribution of the disturbance $\breve{\mathbb{P}}_{w} \in \mathcal{P}^{w}$ are different from their respective worst case counterparts from the assumption of the lemma, it is evident that the state error $e_{k} \neq 0, \forall k = 0, \dots, T$. This means that the bound on $\mathbf{R}_{\alpha}(T, \pi^{\star})$ given by \eqref{eqn_cvar_upper_bound_bigger} grows with the time $T$ and a tighter risk level $\alpha$. 
\end{proof}

\section{Numerical Results} \label{sec_num_sim}
We consider the vehicle steering problem dynamics described in \cite{aastrom2021feedback}:
\begin{align}
A = \begin{bmatrix} 1 & 0.2 \\ 0 & 1 \end{bmatrix}, B = \begin{bmatrix} 0.06 \\ 0.20 \end{bmatrix}, \Sigma_{w} = \begin{bmatrix} 0.10 & 0.03 \\ 0.03 & 0.20 \end{bmatrix}, \mu_{x_{0}} = \begin{bmatrix} -4 \\ 4 \end{bmatrix}, \Sigma_{x_{0}} = \begin{bmatrix} 0.20 & 0.02 \\ 0.02 & 0.20 \end{bmatrix}, 
\end{align}

\begin{figure}[h]
    \centering
    \begin{minipage}{.49\textwidth}
        \centering
        \includegraphics[scale=0.19]{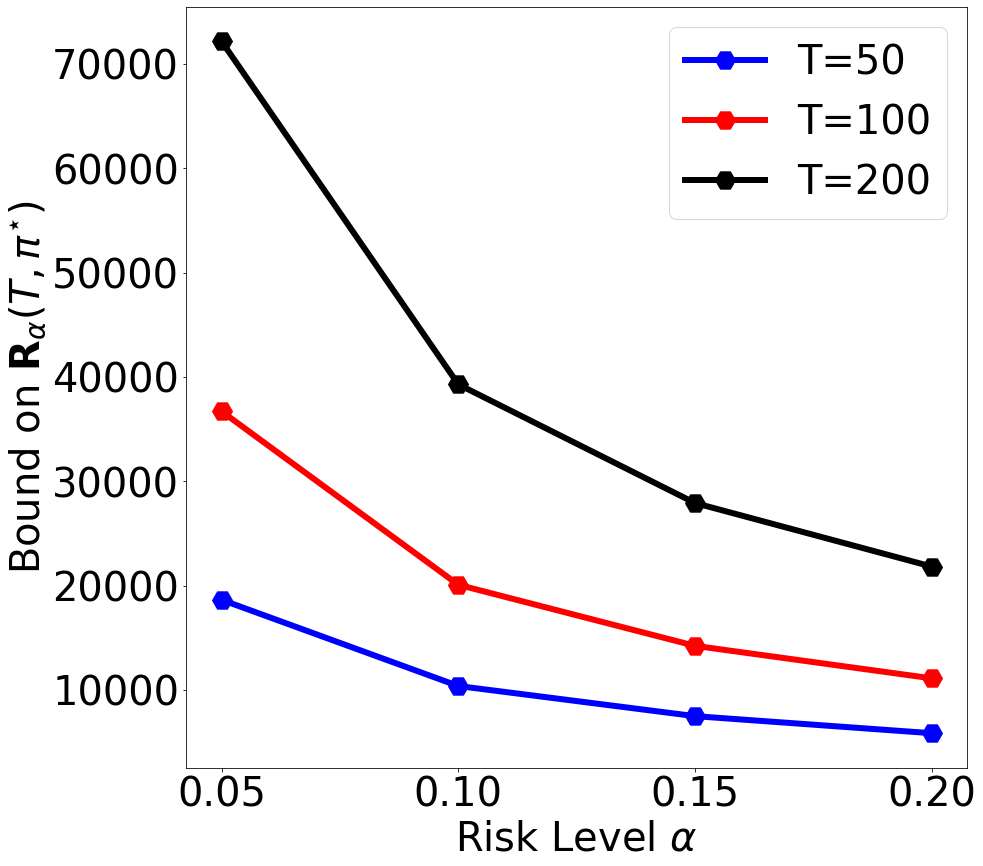}
        \caption{The bounds on the distributional regret $\mathbf{R}_{\alpha}(T, \pi^{\star})$ varying with the risk level $\alpha \in (0,1)$ for different values of time $T$ is shown here. Clearly, the figure depicts that the distributional regret bound defined in  \eqref{eqn_cvar_upper_bound_bigger} increases with time and with tighter risk level.}
        \label{fig_regret_bound}
    \end{minipage}
    \begin{minipage}{.49\textwidth}
        \centering
        \includegraphics[scale=0.19]{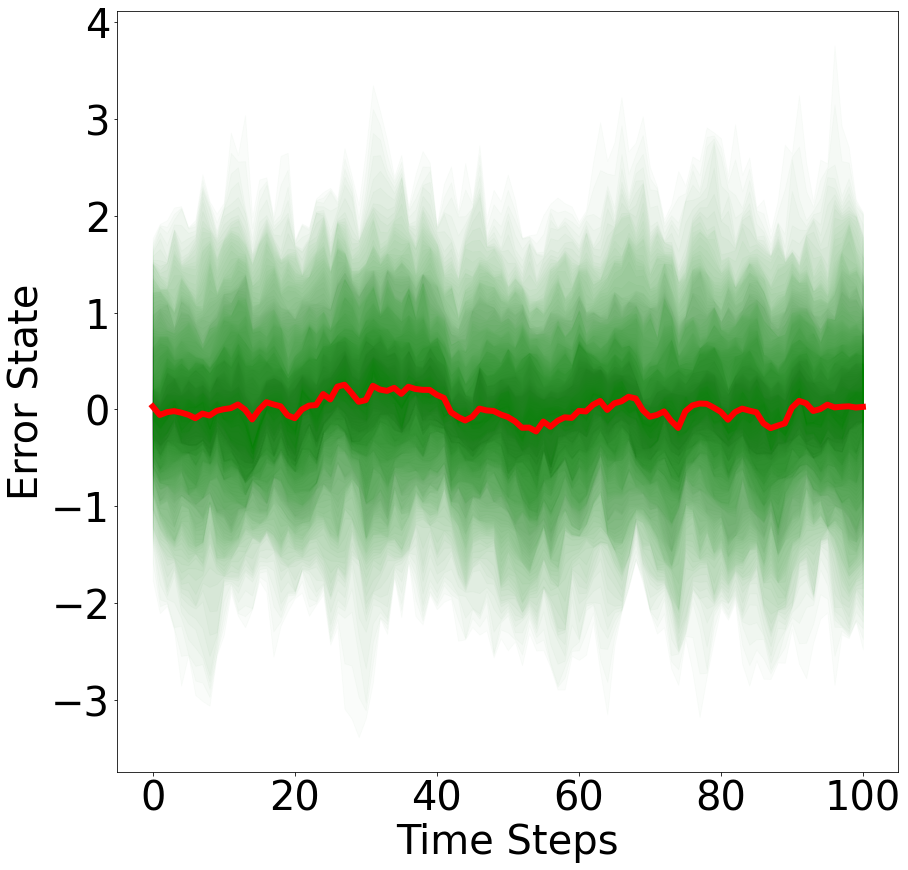}
        \caption{The evolution of the mean of $e_{k}$ and its percentile data are plotted here for time $T = 100$ and $\alpha = 0.20$. The percentile data in shaded green pictorially show the price paid by the policy $\pi^{\star}$ for not being prepared against the uncertainties stemming from the true distributions.}
        \label{fig_error_states}
    \end{minipage}
\end{figure}
\noindent and $D = I_{2}$. Assuming the system state distribution to be stationary in nature, the policy $\pi^{\star}$ was calculated using Theorem 4.2 of \cite{kishida2022risk} where the system was simulated with uncertainties stemming from \emph{possibly} worst-case multivariate Gaussian in the first case and uncertainties stemming from the true distribution namely multivariate Laplacian (as it has heavier tails than the Gaussian counterpart) in the second case. The penalty matrices were assumed to be $Q = 10 I_{2}, R = 1$. The risk level was chosen to be $\alpha = 0.20$. To compute the state error, each time the system was simulated for $T = 100$ time steps using the policy $\pi^{\star}$ from initial state sampled from multivariate Gaussian and multivariate Laplacian distributions respectively. The distributional regret bound corresponding to the term $G_{k}$ was computed using $N = 100$ samples of $e_{k}$. \\

\noindent \emph{Discussion of Results:}
The simulation results are presented in Figures \ref{fig_regret_bound} and  \ref{fig_error_states}. To qualitatively study the bound on the distributional regret defined using Theorem \ref{thm_main}, we simulated the system for different values of the risk level $\alpha \in (0,1)$ by fixing the time $T$. It is clear from the results plotted in Figure \ref{fig_regret_bound} that the distributional regret bound increased with the time $T$ and with a tighter (or smaller) risk level $\alpha$. This means that when the true distribution has heavier tails than the worst-case distribution, the regret incurred will be higher when the risk level $\alpha$ is smaller as the worst-case $\mathrm{CVaR}_{\alpha}(\cdot)$ will increase due to Lemma \ref{lemma_cvar_of_quadratic_cost} and Lemma \ref{lemma_worst_cvar_linear_term}. The evolution of the mean of state error (almost zero for all time steps) and the percentiles are plotted in Figure \ref{fig_error_states} and this depicts that the system evolution with respect to uncertainties from true distribution and worst case distribution is different and the policy that is meant to work for the worst case uncertainty distributions need not work well for the uncertainties stemming from the true distributions governing the uncertainties. The code is available at \url{https://gitlab.control.lth.se/regler/regret-for-unknown-distribution/-/tree/main}.  

\section{Conclusion} \label{sec_conclusion}
The regret incurred by distributionally robust optimal controller while controlling systems with stochastic uncertainties modelled using moment based ambiguity set was presented in this paper. The worst-case $\mathrm{CVaR}_{\alpha}(\cdot)$ with risk level $\alpha \in (0,1)$ was selected to hedge against the distributional uncertainty and using it the regret incurred by the controller for not hedging against the true distribution of system uncertainties was analysed. The distributional regret bound was found to be increasing with tighter (that is smaller) risk levels. Future research will address the distributional regret for risk-averse MPC problem and other DRO parameterizations such as the Wasserstein metric-based ambiguity sets.

\acks{This project has received funding from the European Research Council (ERC) under the European Union’s Horizon 2020 research and innovation program under grant agreement No 834142 (Scalable Control). The authors are with the ELLIIT Strategic Research Area at Lund University. The authors thank Bo Bernhardsson from the department of automatic control in Lund University for his insightful discussions and suggestions.}


\bibliography{sample}

\end{document}